\documentclass{article}
\usepackage{amsmath,amssymb}

\newcommand{\tr}{{\mathrm{tr\,}}}
\newcommand{\FS}{{\mathrm{FS}}}
\newcommand{\Sm}{{\mathrm{Sm}}}
\newcommand{\Br}{{\mathrm{Br}}}
\newcommand{\fix}{{\mathrm{fix}}}
\newcommand{\Smo}{{\preceq_{\Sm}}}
\newcommand{\tSmo}{{\preceq_{\underline\Sm}}}
\newcommand{\Bro}{{\preceq_{\Br}}}
\newcommand{\Smr}{{\leftarrow_{\Sm}\,}}
\newcommand{\tSmr}{{\leftarrow_{\underline \Sm}\,}}

\newcommand{\R}{{\mathbb R}}
\newcommand{\C}{{\mathbb C}}

\newcommand{\Proj}{{\mathbb P}}
\newcommand{\ad}{{\rm ad}}
\newcommand{\spec}{{\rm spec}}

\newcommand{\GL}{{\rm GL}}

\newcommand{\gl}{{\mathfrak g \mathfrak l}}

\newcommand{\g}{{\mathfrak g}}

\newcommand{\kk}{{\mathfrak k}}
\newcommand{\pp}{{\mathfrak p}}
\newcommand{\aaa}{{\mathfrak a}}
\newcommand{\m}{{\mathfrak m}}
\newcommand{\n}{{\mathfrak n}}

\newcommand{\bb}{{\mathfrak b}}

\newtheorem{PARA}{\S}

\newcommand{\MAT}[1]{\left[\begin{array}{#1}}
\newcommand{\RIX}{\end{array}\right]}
\newcommand{\End}{\rm End}
\newcommand{\la}{\langle}
\newcommand{\ra}{\rangle}

\newtheorem{theorem}{Theorem}[section]
\newtheorem{remark}[theorem]{Remark}
\newtheorem{lemma}[theorem]{Lemma}
\newtheorem{corollary}[theorem]{Corollary}
\newtheorem{definition}[theorem]{Definition}
\newtheorem{proposition}[theorem]{Proposition}

\newcommand{\halmos}{\rule{1.5mm}{1.5mm}}
\newenvironment{proof}[0]{\noindent {\it Proof.\small}
}{\hfill \halmos\vskip 1cm
}

\begin{document}

\title{
Attractor Networks on Complex Flag Manifolds}
\author{
Joachim Hilgert\\ Institut f\"ur Mathematik\\ Universit\"at
Paderborn\\ D-33095 Paderborn\\ Germany }
\date{}
\maketitle

\begin{abstract} Robbin and Salamon showed in \cite{RobbinSalamon}
that attractor-repellor networks and Lyapunov maps are equivalent
concepts and illustrate this with the example of linear flows on
projective spaces. In these examples the fixed points are linearly
ordered with respect to the Smale order which makes the
attractor-repellor network overly simple. In this paper we provide a
class of examples in which the attractor-repellor network and its
lattice structure can be explicitly determined even though the Smale
order is not total. They are associated with special flows on
complex flag manifolds. In the process we show that the Smale order
on the set of fixed points can be identified with the well-known
Bruhat order. This could also be derived from results of Kazhdan and
Lusztig, but we give a new proof using the $\lambda$-Lemma of Palis.
For the convenience of the reader we also introduce the flag
manifolds via elementary dynamical systems using only a minimum of
Lie theory.
\end{abstract}

\pagenumbering{arabic}

\section{Introduction}\label{introduction}

Given a dynamical system  on a space ${\cal M}$, Smale
considered the relation $x\Smr y$ defined by
${\cal W}^-(y)\cap {\cal W}^+(x)\not=\emptyset$, where ${\cal W}^\pm(m)$ are the
stable and unstable manifolds of $m\in{\cal M}$. This relation in general is
not transitive, so  one takes the transitive closure which is then called
the Smale order and denoted by $\Smo$.

Let now ${\cal M}$ be a complex flag manifold. Viewing ${\cal M}$ as
an adjoint orbit it is possible to construct gradient flows on
${\cal M}$ such that the Bruhat cells in ${\cal M}$ coincide with
the unstable manifolds of this flow as Atiyah remarks in
\cite{Atiyah}. The set of Bruhat cells carries a natural order
defined by $N_1\Bro N_2$ if $N_1\subseteq \overline{N_2}$, where
$\overline{N_2}$ is the closure of $N_2$ in ${\cal M}$. This order
is called the Bruhat order and it actually coincides with the
combinatorial Bruhat order on the coset space of the Weyl group
associated with ${\cal M}$  (see \cite[Lemma~6.2.1]{BE}). This
observation suggests a close relation between the Smale order and
the Bruhat order. In fact, it turns out that the relation $\Smr$ in
this case is a partial order which coincides with $\Bro$ when the
Weyl group coset space is identified with the set of fixed points.
This allows us to describe the algebraic structure of the attractor
network associated in \cite{RobbinSalamon} with the dynamical system
on ${\cal M}$ in terms of the well studied order structures on Weyl
groups.

The proofs we present are based on linear flows on projective space. The elementary techniques involved can
also be used to provide quick proofs for most of the basic results on flag manifolds required in our context.
We included them for the convenience of the reader not specializing in Lie theory.

\subsection*{Acknowledgement}

The author would like to thank J.W.~Robbin for bringing attractor
networks and Lyapunov maps to his attention. The ensuing discussions
of the projective space examples triggered the idea to construct
more complicated examples using flag manifolds.

\section{The Smale Order}

Let ${\cal M}$ be  a topological space and $\Phi\colon \R\times
{\cal M}\to {\cal M}$ a continuous group action  of $\R$ on   ${\cal
M}$. Then we call $\Phi$ a dynamical system with continuous time. We
assume that ${\cal M}$ is compact and admits a metric $d$. For any
invariant subset ${\cal A}\subseteq {\cal M}$  the {\it stable} and
{\it unstable manifold} of ${\cal A}$ in ${\cal M}$ are defined by
$${\cal W}^\pm({\cal A}):=
\{x\in {\cal M}\vert \lim_{t\to\pm \infty}d(\Phi(t,x),{\cal A})=0\}.$$
Note that ${\cal W}^\pm({\cal A})$
are both invariant under the flow.
One defines a relation $\Smr$ on the set
of fixed points ${\cal M}_\fix$ by
\[x\Smr y\quad:\Leftrightarrow\quad
{\cal W}^-(y)\cap {\cal W}^+(x)\not=\emptyset.\]
Then the  transitive closure $\Smo$
of $\Smr$ is called the {\it Smale order}.
Note that  in general the relation $\Smo$
is {\it not} antisymmetric, i.e. it is {\it not} a partial order.
In order to be able to derive a reasonable partial order from
$\Smo$ one has to make some additional assumptions
on the dynamical system.

\begin{definition}\label{admissible}{\rm
We call a dynamical system $({\cal M},\Phi)$ {\it admissible} if
${\cal M}$ is compact, every flow line has a sink and a source,
and, in addition, ${\cal M}_{\fix}$ has only finitely many open
closed subsets,
i.e.\ ${\cal M}$ has only finitely many compact isolated
subsets of fixed points. }
\end{definition}

Now suppose that $({\cal M},\Phi)$ is admissible. We denote the
(finite) set of connected open closed subsets of ${\cal M}_{\fix}$
by $\underline{\cal M}_{\fix}$. Then we can define a {\it Smale
order} $\tSmo$ also on $\underline{\cal M}_{\fix}$  as the
transitive closure of the relation
\[p\tSmr q\quad:\Leftrightarrow\quad
{\cal W}^-(q)\cap {\cal W}^+(p)\not=\emptyset\] on $\underline{\cal
M}_{\fix}$.

\begin{remark}\label{tildeM3}{\rm The set  of subsets
$F\subseteq \underline {\cal M}_{\fix}$ is in bijective
correspondence with the set of open closed subsets ${\cal
U}\subseteq{\cal M}_{\fix}$ via
$$F\mapsto {\cal U}_F:= \bigcup_{p\in F}p,\quad
{\cal U}\mapsto F_{\cal U}:=
          \{p\in\underline{\cal M}_{\fix}\mid p\subseteq {\cal U}\}.$$
A set $F\subseteq \underline{\cal M}_{\fix}$ is an upper set with
respect to the Smale order if and only if
$${\cal W}^+({\cal U}_F)\cap {\cal M}_{\fix}={\cal U}_F$$
which means that ${\cal U}_F$ is an upper set in ${\cal M}_{\fix}$.
Similarly, $F\subseteq \underline{\cal M}_{\fix}$ is a lower set
with respect to the Smale order if and only if  ${\cal U}_F$ is a
lower set in ${\cal M}_{\fix}$.}
\end{remark}

\section{Linear Flows on Projective Space}
\label{LinFlow}

Let $V$ be finite dimensional complex vector space with an inner
product $\la\cdot\,|\cdot \ra$ and
 $\varphi\in \End(V)$ a selfadjoint linear map.
We denote the (real) spectrum of $\varphi$ by $\spec(\varphi)$ and
the eigenspace of $\varphi$ for $\lambda\in \spec(\varphi)$ by
$V_\lambda$. We consider the flow induced by $e^{t\varphi}$ on the
projective space $\Proj(V)$ the  projective flow generated by
$\varphi$  and denote it by $\Phi$.

\label{gradient} The  flow $\Phi$ of $\varphi$ on $\Proj(V)$ can be
interpreted as a gradient flow. In fact, denote the natural
projection $V\setminus\{0\}\to \Proj(V)$ by $\pi$ and equip
$\Proj(V)$ with the  Fubini-Study metric $\la\cdot\,|\cdot
\ra_{\FS}$, which turns $\Proj(V)$ into a K\"ahler manifold and
which is characterised by the equation
\[\la d\pi_v(u)\mid d\pi_v(w)\ra_{FS}=
\frac{\la u\mid w\ra}{\la v\mid v\ra}\]
at the point $[v]$, where $u$ and $w$ are orthogonal to $v$. Then we have

\begin{proposition} \label{gradient1}
 $[v]\mapsto -d\pi_v(\varphi v)$ is the gradient vector
field of the height function
\[f([v])=-\frac{1}{2}\frac{(\varphi v\mid v)}{(v\mid v)},\]
where $(v\mid w)=\Re\la v\mid w\ra$ and the Riemannian metric on
$\Proj(V)$ is the real part of the Fubini-Study metric, i.e.,
\[\nabla f([v])=-d\pi_v(\varphi v).\]
\end{proposition}

\begin{proof} An elementary calculation shows that
\[d(f\circ \pi)_v(w)=\frac{-1}{(v\mid v)}
\Big(\varphi v-\frac{(\varphi v\mid v)}{(v\mid v)}v\ \Big|\  w\Big).\]
On the other hand we have
\[d(f\circ \pi)_v(w)=\Re\la \nabla f([v])\mid d\pi_v(w)\ra_{\FS}.\]
Since $v$ is in the kernel of $d\pi_v$ the claim follows.
\end{proof}

Note that Proposition \ref{gradient1} in particular shows that $f$
is strictly decreasing along non-constant flow lines.
\label{fixedpoints} The fixed points of $e^{t\varphi}$ on $\Proj(V)$
are
\[{\cal F}:=\Proj(V)_{\fix}=\{[v]\in \Proj(V)\mid
\big(\exists \nu\in \spec(\varphi)\big)\, v\in V_\nu\},\] where
$[v]$ is the line spanned by $v$. There is a natural embedding of
$\Proj(V_\nu)\to \Proj(V)$ for each $\nu\in \spec(\varphi)$. The
height function takes the value $f([v])=-\frac{1}{2}\nu$ on
$\Proj(V_\nu)$. Identifying $\Proj(V_\nu)$ with its image under this
embedding we obtain
\begin{equation}\label{tildeM1}
{\cal F}=\bigcup_{\nu\in \spec(\varphi)} \Proj(V_\nu).
\end{equation}
Given $x=[v]\in \Proj(V)$ we define
\[w_+(x):={\rm min}\{r\in \R\mid v\in
\bigoplus_{\nu\le r} V_\nu\}, \quad w_-(x):={\rm max}\{r\in \R\mid
v\in \bigoplus_{r\le \nu} V_\nu\} .\] In particular, we have
$w_\pm(x)=\nu$ for $x\in \Proj(V_\nu)$. In this case we simply write
$w(x)$ for $w_+(x)$ and call it the $\varphi$-{\it weight} of $x$.

\begin{proposition}\label{sinksourceproj}
Each $e^{\R\varphi}$--orbit has  limits for $t\to\pm\infty$. More
precisely, for $x=[v]$ with $v=\sum_{\nu\in \spec(\varphi)}v_\nu$
and $v_\nu\in V_\nu$ we have $\lim_{t\to \infty}e^{\pm
t\varphi}\cdot x=[v_{w_\pm(x)}]$. We denote the resulting map
$\Proj(V)\to {\cal F}\times {\cal F}$ by
\[\ell\colon \Proj(V)\to {\cal F}\times {\cal F},\quad
x\mapsto
\lim_{t\to\infty}(e^{-t\varphi}\cdot x,e^{t\varphi}\cdot x).\]
\end{proposition}

\begin{proof} This is a straightforward calculation in
homogeneous coordinates with respect to any basis which diagonalises
$\varphi$.
\end{proof}

\begin{remark}\label{tildeM2}{\rm
Proposition \ref{sinksourceproj} together with (\ref{tildeM1}) shows
that $(\Proj(V),\Phi)$ is admissible and
$\underline{\Proj(V)}_{\fix}$ with the Smale order is order
anti-isomorphic to $\spec(\varphi)$ with the natural order.}
\end{remark}

\label{stableunstable}
It is an elementary exercise to calculate
the stable manifold and unstable manifolds:

\begin{proposition}
Let $x=[v_\nu]\in \Proj(V_\nu)\subseteq {\cal F}$ for some $\nu\in
\spec(\varphi)$. Then
\[{\cal W}^\pm(x)=
\pi\left(v_\nu+\sum_{0<\pm(\nu-\mu)}V_\mu\right)
.\]
\end{proposition}

\section{Complex Flag Manifolds}
\label{flagmanifolds}

In this section we review the construction and some basic properties of
complex flag manifolds. The results presented are standard but the
approach via dynamical systems is elementary and not so widely known.
Note, however, that Duistermaat, Kolk, and Varadarajan in \cite{DKV}
provide plenty of information on this approach also in the more complicated
case of real flag manifolds.

\label{lingroups}
Let  $V$ be a finite dimensional complex vector space and
$G\subseteq \GL(V)$ be a closed  connected
subgroup. The {\it Lie algebra} of $G$ is the set
\[\g:=\{X\in {\End}_\C(V)\mid e^{\R X}\subseteq G\}.\]
Then $\g$ is a real vector subspace of
$\gl(V):=\End_\C(V)$ which is closed under the bilinear
map
\[\gl(V)\times \gl(V)\to\gl(V),\quad (X,Y)\mapsto
[X,Y]:=XY-YX\]
called the {\it Lie bracket}. The matrix exponential map
$\exp:\g\to G, \ X\mapsto e^X$ is a
local homeomorphism at $0$ and can be used to define
a $G$-invariant differential structure on $G$ such that
$G$ is a closed real submanifold of the open subset
$\GL(V)$ in $\End_\C(V)$.
If  $\g$ is a complex vector subspace of $\gl(V)$,
then $G$ is a complex submanifold of $\GL(V)$.
The group generated by $\exp(\g)$ is all of $G$.

Suppose now that $V$ is equipped with an inner
product  $(v,w)\mapsto \la v\mid w\ra$.
Given $\varphi\in \End_\C(V)$ let $\varphi^t\in \End_\C(V)$
be the transpose of $\varphi$ with respect to the inner product.
The map
\[\theta:\GL(V)\to\GL(V),\quad g\mapsto (g^{-1})^t\]
is involutive and preserves the group multiplication.
It is called the {\it Cartan involution}. Its derivative at the
identity is the map
\[ \gl(V)\to\gl(V),\quad X\mapsto -X^t\]
which  preserves the Lie bracket. By abuse of notation
it is also called Cartan involution and denoted by $\theta$.
If $G$ is $\theta$--invariant, then it
is easy to check that $\g$ is invariant under
the (algebra) Cartan involution, so that we have a
{\it Cartan} decomposition
$\g=\kk+\pp$,
where $\kk$ and $\pp$ are the $\theta$--eigenspaces
for the eigenvalues $1$ and $-1$, respectively.  The polar decomposition yields a global
analogue, also called Cartan decomposition:
$G=KP$. Here $K$ is the subgroup of $\theta$--fixed points in $G$
and $P=\exp(\pp)$.

\label{weightsandroots}
A subspace of $\gl(V)$ is called {\it abelian} if the Lie
bracket vanishes on it. Choose a maximal abelian
subspace $\aaa$ in $\pp$. Then $\aaa$ consists of
commuting $\R$--split endomorphisms and is therefore
simultaneously diagonalizable. Thus there exists a set
${\cal P}$ of linear functionals on $\aaa$ such that
\[V=\bigoplus_{\mu\in {\cal P}}V^\mu,\]
where
$V^\mu=\{v\in V\mid (\forall X\in \aaa)\,  Xv=\mu(X)v\}$
and $V^\mu\not=\{0\}$ for all $\mu\in {\cal P}$.
The elements of ${\cal P}$ are called the {\it weights}
of $G$ on $V$.
Each element $X$ of $\g$ defines a linear map
$\ad(X):\g\to \g$ via
$\ad(X)Y=[X,Y]$, which is called the {\it adjoint representation} of $X$.
Note that $\ad(X)$ is self adjoint with respect to the inner product
$(Y,Z)\mapsto \tr(YZ^t)$ on $\g$ whenever
$X\in \pp$. Therefore these $\ad(X)$ are $\R$--split.
Since $\aaa$ is abelian the maps $\ad(X)$ with
$X\in \aaa$ commute. Thus, in analogy to the weight
decomposition of $V$, we have a {\it root decomposition}
$$\g=\m\oplus\aaa\oplus\bigoplus_{\alpha\in \Delta}
\g^\alpha,$$
where
$\g^\alpha:=\{Y\in \g\mid (\forall X\in \aaa)\ \ad(X)Y=\alpha(X)Y\}$,
$\Delta:=\{\alpha\in \aaa^*\mid {\alpha\not=0,} \ {\g^\alpha\not=\{0\}}\}$,
and
$\m:=\{Y\in \kk\mid (\forall X\in \aaa)\ [X,Y]=0\}$.
The elements of $\Delta$ are called the
{\it restricted roots}
associated with the pair $(\g,\aaa)$. Note that
$\theta(X)=-X$ for $X\in \aaa$
implies that for $\alpha\in \Delta$ also $-\alpha$ is a
restricted root
and $\theta\g^\alpha=\g^{-\alpha}$.
Choose an element
$X_o\in \aaa$ such that $\alpha(X_o)\not=0$ for all
$\alpha\in \Delta$
and set
$\Delta^\pm:=\{\alpha\in \Delta\mid \pm\alpha(X_o)>0\}$.
the elements of $\Delta^+$ are called  {\it positive roots}
and  $\Delta^-=-\Delta^+$ is the set of {\it negative roots}.
The Lie bracket satisfies the {\it Jacobi identity}
$ [X,[Y,Z]]+[Y,[Z,X]]+[Z,[X,Y]]=0$
from which one derives
\[\g^\alpha(V^\mu)\subseteq V^{\mu+\alpha},
\quad [\g^\alpha,\g^\beta]\subseteq \g^{\alpha+\beta}.\]
Therefore
$\n:=\bigoplus_{\alpha\in \Delta^+}\g^\alpha$
 is a subalgebra of $\g$ consisting of
simultaneously triangularizable elements of $\End_\C(V)$
and we obtain the {\it Bruhat} decomposition
\[\g=\theta\n\oplus\m\oplus \aaa\oplus\n.\]
On the other hand $Y+\theta(Y)\in \kk$ for $Y\in \theta\n$
so that $\theta\n\subseteq \kk+\n$. This yields
the {\it Iwasawa} decomposition
\[\g=\kk+\aaa+\n\]
which is easily checked to be direct using the
properties of $\theta$.
 The
global version of the Iwasawa decomposition is
$G=KAN$, where $A=\exp(\aaa)$ and $N=\exp(\n)$.

\label{Borel}
The {\it commutator algebra} of a Lie algebra $\g$ is the space $\g'$
spanned by all $[X,Y]$ with $X,Y\in \g$. It is easy to check that it is actually a
subalgebra of $\g$ under the Lie bracket. One can show
that the Lie algebra of the commutator subgroup $G'$ of
$G$ is $\g'$. The Lie
algebra $\g$ is called {\it solvable} if the sequence of
successive commutator algebras
\[\g\supseteq \g'\supseteq (\g')'\supseteq\ldots\]
reaches zero.

Suppose now that $\g$ is a complex subspace of $\gl(V)$. Then
$\kk=i\pp$ and hence $\aaa_\C=\aaa+i\aaa$ is maximal abelian in
$\g$. As a consequence $\bb:=\m+\aaa+\n$ is a complex solvable
subalgebra of $\g$ called the {\it Borel} subalgebra. The (closed)
subgroup of $G$ generated by $\exp\bb$ is called {\it Borel
subgroup} and denoted by $B$. The group $\exp(\aaa_\C)$ is denoted
by $A_\C$. The abelian algebra $\aaa_\C$ is its own normaliser in
$\g$, i.e., $[X,\aaa_\C]\subseteq \aaa_\C$ implies $X\in \aaa_\C$.
Such algebras are called {\it Cartan subalgebras} and one can show
that all Cartan subalgebras of $\g$  are conjugate under $G$. Using
Jordan canonical forms, linear flows on $\Proj(V)$  can be used to
give an elementary proof of the following version of Borel's Fixed
Point Theorem:

\begin{theorem} \label{Borel1}
Let $S\subseteq \GL(V)$ be a closed connected subgroup with complex
solvable Lie algebra  and $C\subseteq \Proj(V)$ a closed
$S$--invariant subset. Then $C$ contains an $S$-fixed point.
\end{theorem}

Let ${\cal M}:=G\cdot x\subseteq \Proj(V)$ be a closed $G$--orbit.
Then the above theorem applied to the group $B$
 shows that there exists a $B$--fixed point
$x_c\in {\cal M}$. In other words, $B$ is contained in the
stabiliser $G_{x_c}$ of $x_c$.

\label{parabolics} Subgroups of $G$ containing $B$ are called {\it
standard parabolic} subgroups of $G$. Conjugates of standard
parabolic subgroups are called {\it parabolic} subgroups. Thus all
the stabilisers $G_y$ of points $y\in {\cal M}$ are parabolic
subgroups. Conversely, suppose that $B\subseteq G_y$ for some $y\in
\Proj(V)$.
 The Iwasawa decomposition shows that
$G/G_y\cong K/K_y$ is compact so that $G\cdot y$ is closed in
$\Proj(V)$. One can even show that all parabolic subgroups can be
obtained in this way. A {\it complex flag manifold} of $G$  is a
homogeneous space of the form $G/P$ with $P$ parabolic. Thus our
discussion shows that the complex flag manifolds are precisely the
closed $G$-orbits in $\Proj(V)$.

We close this section with a technical  lemma which  is the  tool we
need to characterise the fixed point set of our flows on complex
flag manifolds. \label{Weylgroup} Let $N_K(\aaa)=\{k\in K\mid k\aaa
k^{-1}=\aaa\}$ and $Z_K(\aaa)=\{k\in K\mid (\forall X\in \aaa)\
kXk^{-1}=X\}$ be the normaliser and the centraliser of $\aaa$ in
$K$, respectively.

\begin{lemma} \label{Weylgroup1}
\begin{enumerate}
\item[{\rm (i)}] $Z_K(\aaa)$ is a normal subgroup of
$N_K(\aaa)$ and the  quotient group
$W:=N_K(\aaa)/Z_K(\aaa)$ is finite.
\item[{\rm (ii)}] Suppose that $x,y\in {\cal M}$ are both
fixed under $A_\C$. Then there exists a $k\in N_K(\aaa)$
such that $k\cdot x=y$.
\end{enumerate}
\end{lemma}

\begin{proof} (i) The first part is obvious. To prove the second
part note first that the annihilator $\Delta^\perp$ of $\Delta$ in $\aaa$ is
central in $\g$. The group $N_K(\aaa)$ permutes the elements
of $\Delta$ (dual action). If an element $k\in N_K(\aaa)$
fixes all $\alpha\in \Delta$, then $k$ acts trivially on
$\aaa/\Delta^\perp$ and by the above remark also on $\aaa$.

(ii) Using the global Iwasawa decomposition again, we see that $K$
acts transitively on ${\cal M}$, so we find a $k_1\in K$ such that
$y=k_1\cdot x$. Therefore the stabilisers of $x$ and $y$ are
conjugate under $k_1$, i.e. $G_y=G_{k_1\cdot x}=k_1G_x k_1^{-1}$.
This shows that $A_\C$ and $k_1^{-1}A_\C k_1^{-1}$ are contained in
$G_x$. Therefore $\aaa_\C$ and $k_1\aaa_\C k_1^{-1}$ are contained
in the Lie algebra $\g_y$  of $G_y$ and hence automatically Cartan
subalgebras of $\g_y$. Now we use the fact that all Cartan
subalgebras of $\g_y$ (which is complex) are conjugate under $G_y$.
This means we can find an $h\in G_y$ such that
$\aaa_\C=hk_1\aaa_\C(hk_1)^{-1}$. Then $hk_1$ is contained in the
$\theta$--invariant set $N_G(\aaa_\C)$, which according to the
global Cartan decomposition, is equal to $N_K(\aaa)A$. Therefore we
can write $hk_1=ka_1$ with $a_1\in A$ and obtain
\[k\aaa k^{-1}=hk_1a_1^{-1}\aaa a_1(hk_1)^{-1}=hk_1\aaa (hk_1)^{-1}=\aaa\]
as well as
$k\cdot x=hk_1a_1\cdot x=hk_1\cdot x=h\cdot y=y$.
\end{proof}

The group $W=N_K(\aaa)/Z_K(\aaa)$ is called the {\it Weyl group} of
$G$. One can show that $Z_K(\aaa)$ is connected (being the
centraliser in $K$ of the maximal abelian subalgebra $i\aaa$ in
$\kk$, cf.\ \cite[p.~287]{Helgason78}). Therefore it is equal to
$\exp(i\aaa)$ and thus contained in $B$. Thus we can define a
subgroup $W_{\cal M}$ of $W$ via
\[W_{\cal M}:=\{kZ_K(\aaa)\mid k\cdot x_c=x_c\}.\]

\section{Special Flows on Complex Flag Manifolds}
\label{specialflows}

In this section we introduce the kind of flow we want to
study on complex flag manifolds.

\begin{definition}\label{specialflow}
Suppose that $X\in \aaa$ has the following properties:
\begin{enumerate}
\item[{\rm (i)}] $\alpha(X)< 0$ for all $\alpha\in \Delta^+$.
\item[{\rm (ii)}] $\mu(X)\not=\nu(X)$ for all $\mu,\nu\in {\cal P}$
with $\mu\not=\nu$.
\end{enumerate}
Then we call the  linear flow $\Phi_X$ of $X$ on $\Proj(V)$ {\it
special}.
\end{definition}

Suppose that the flow of $X\in \aaa$ is special. Then the  fixed
points of $\Phi_X(\R)=e^{\R X}$ in $\Proj(V)$ are precisely the
points $[v]$ with $v\in V^\mu$ for some $\mu\in {\cal P}$. In fact,
$[v]$ is a fixed point for $e^{\R X}$ if and only if $\C v$ is
invariant under $X$ which in turn is the same as saying that $v$ is
an eigenvector of $X$. But our hypothesis on $X$ implies that the
eigenvalues of $X$ are the $\mu(X)$ with $\mu\in {\cal P}$ and that
$\mu(X)$ determines $\mu$. Thus any eigenvector of $X$ belongs to
some $V^\mu$ with $\mu\in {\cal P}$. As a consequence we see that an
$e^{\R X}$--fixed point in $\Proj(V)$ is automatically an
$A_\C$-fixed point, where $A_\C=\exp(\aaa_\C)$. \label{Bruhatcells}
We restrict the special flow $e^{\R X}$ to ${\cal M}=G\cdot
x_c\subseteq \Proj(V)$ and consider the set
\[{\cal F}_{\cal M}={\cal F}\cap {\cal M}
=\{x\in {\cal M}\mid A_\C\cdot x=x\}\]
of $A_\C$--fixed points in ${\cal M}$ as well as the
corresponding stable and unstable manifolds
${\cal W}_{\cal M}^{\pm}(x)=
{\cal W}^{\pm}(x)\cap {\cal M}$
in ${\cal M}$.

\begin{theorem} \label{Bruhatcells1}
\begin{enumerate}
\item[{\rm (i)}] The Weyl group $W$ acts
transitively  on ${\cal F}_{\cal M}$. The stabiliser of $x_c$ in $W$
is $W_{\cal M}$. Therefore ${\cal F}_{\cal M}$ is parametrised by
$W/W_{\cal M}$ and $({\cal M},\Phi)$ is admissible with
$\underline{\cal M}_{\fix}={\cal M}_{\fix}={\cal F}_{\cal M}$.
\item[{\rm (ii)}] ${\cal M}=
\bigcup_{x\in {\cal F}_{\cal M}} {\cal W}^+(x)=
\bigcup_{x\in {\cal F}_{\cal M}} {\cal W}^-(x)$.
\item[{\rm (iii)}] ${\cal W}_{\cal M}^+(x)=N\cdot x$ and
${\cal W}_{\cal M}^-(x)=(\theta N)\cdot x$
for all $x\in {\cal F}_{\cal M}$.
\end{enumerate}
\end{theorem}

\begin{proof} (i) The first part follows from   Lemma \ref{Weylgroup1}.
The second part is an immediate
consequence of the definition of $W_{\cal M}$.

(ii) Since all eigenvalues of $X$ are real, the
limit of $e^{tX}\cdot y$ for $t\to\pm\infty$
exists and is contained in ${\cal F}_{\cal M}$ by
Proposition \ref{sinksourceproj}.

(iii) Recall that
$N=\exp(\sum_{\alpha\in \Delta^+}\g^\alpha)$.
We calculate
\begin{eqnarray*}
\exp(tX)\exp(\sum Y^\alpha)\cdot x&=&
\exp(tX)\exp(\sum Y^\alpha)\exp(-tX)\cdot x\\
&=&\exp(\sum e^{\ad tX}Y^\alpha)\cdot x\\
&=&\exp(\sum e^{t\alpha(X)}Y^\alpha).x
\end{eqnarray*}
and note that this expression converges to
$x$ for $t\to \infty$ by our hypothesis on $X$. As a result
we obtain
$N\cdot x\subseteq {\cal W}_{\cal M}^+(x)$.
Analogously we find
$\theta N\cdot x\subseteq {\cal W}_{\cal M}^-(x)$.

To show the converse, note first that the tangent space
$T_x{\cal M}$ is given by
\[T_x{\cal M}=\{X.x\mid X\in \theta\n+\n\}\]
since $A_\C$ fixes $x$. Therefore the above calculation
shows that locally (close to $x$) $N\cdot x$ agrees with
${\cal W}_{\cal M}^+(x)$
and $\theta N\cdot x$ agrees with ${\cal W}_{\cal M}^-(x)$.
 Since $\exp(tX)$ normalises $N$ and $\theta N$ this
proves the claim.
\end{proof}

As a corollary we obtain a generalization of the global version of the Bruhat decomposition
which explains why we call the
$N$-orbits in ${\cal M}$  {\it Bruhat cells}:

\begin{corollary}\label{Bruhatdecomp}
Let $m_1,\ldots, m_l\in N_K(\aaa)$ be such that
$$m_1Z_K(\aaa),\ldots , m_lZ_K(\aaa)\in W$$
is a system of
representatives for $W/W_{\cal M}$. Then we have a decomposition
\[G=\bigcup_{j=1}^l Nm_jB.\]
\end{corollary}

The following lemma shows that special flows have desirable properties from
the point of view of dynamical systems. It will allow us to apply the
$\lambda$-Lemma of Palis (see \cite[\S 2.7]{PM}),
which will eventually show the transitivity of the
Smale relation.

\begin{lemma} \label{specialisMorseSmale}
The vector field $X$ on ${\cal M}$ is Morse-Smale.
In fact, we even have:
\begin{enumerate}
\item[{\rm (i)}] $-Xv$ is the gradient vector field of the
function
\[f([v])=-\frac{1}{2}\frac{(Xv\mid v)}{(v\mid v)},\]
where $(v\mid w)=\Re\la v\mid w\ra$ and the Riemannian metric on
$\Proj(V)$ is the Fubini-Study metric.
\item[{\rm (ii)}] ${\cal W}_{\cal M}^-(x)$ and
${\cal W}_{\cal M}^+(y)$ intersect transversally for all
$x,y\in {\cal F}_{\cal M}$.
\item[{\rm (iii)}] All critical points of $X$ are hyperbolic.
\end{enumerate}
\end{lemma}

\begin{proof} (i)  follows immediately from
Proposition \ref{gradient1}.

(ii) All $N$-- and $\theta N$--orbits  of elements in ${\cal
F}_{\cal M}$ are $A_\C$--invariant since $A_\C$ normalises $N$ and
$\Theta N$. Therefore we have
\[T_z{\cal M}=(\n+\theta \n)\cdot z=
T_z(N\cdot z)+T_z(\theta N\cdot z)=
T_z({\cal W}_y^+)+T_z({\cal W}_x^-)\]
for each point $z\in {\cal W}_{\cal M}^-(x)
\cap {\cal W}_{\cal M}^+(y)$.

(iii) This follows from the fact that $X$ is semisimple
with real eigenvalues.
\end{proof}

\section{The Bruhat Order and its Relation to the Smale Order}

\label{Bruhatorder}

Consider the relation
\[x\Bro y\quad :\Longleftrightarrow\quad N\cdot x
\subseteq \overline{N\cdot y}\] on the set ${\cal F}_{\cal M}$ of
$A_\C$--fixed points in ${\cal M}$. It clearly is reflexive and
transitive. We will see below that it actually also is
antisymmetric, i.e., a partial order. This order is called the {\it
Bruhat order} on ${\cal F}_{\cal M}$. Recall from Theorem
\ref{Bruhatcells1} that ${\cal F}_{\cal M}$ is parametrised by
$W/W_{\cal M}$. A result of Kostant (see
\cite[Lemma~1.1.2.15]{Warner72}) shows that there is a canonical
choice of representatives $W^{\cal M}$ for $W/W_{\cal M}$ in $W$ and
using the generalised Bruhat decomposition of $G$ (see
\cite[Chap.~4,\S~2(5), Prop.~2]{Bou81}) one can show that the
induced order on $W^{\cal M}$ actually is the Bruhat order as
defined by generators and relations for $W$ viewed as a Coxeter
group (see \cite[\S~5.10]{Humphreys}  and also
\cite[Lemma~6.2.1]{BE}).

The Smale relation and the Smale order in this context are given by
\begin{eqnarray*}
x\Smr y\quad &\Longleftrightarrow&\quad (\theta N)\cdot y\ \cap \
                                        N\cdot x\not=\emptyset,\\
x\Smo y\quad&\Longleftrightarrow&\quad x\Smr\ldots\Smr y.
\end{eqnarray*}

In view of Lemma \ref{specialisMorseSmale} the $\lambda$-Lemma of
Palis now shows that given $x\in {\cal F}_{\cal M}$ and two
submanifolds $S_{\pm}$ of ${\cal M}$ with $S_{\mp}$ intersecting
${\cal W}_{\cal M}^{\pm}(x)$ transversally in $x_{\mp}$ one can find
a $t>0$ such that $e^{tX}\cdot S_-\cap S_+\not=\emptyset$. This is
the key observation in the proof of our main tool for the
determination of the attractor networks of complex flag manifolds:

\begin{theorem} \label{Bruhatorder1}
\begin{enumerate}
\item[{\rm (i)}] The relations $\Bro$,$\Smr$, and $\Smo$
on ${\cal F}_{\cal M}$ agree.
\item[{\rm (ii)}] $\overline{N\cdot x}=
\bigcup_{y\Bro x} N\cdot y$.
\item[{\rm (iii)}] $\overline{\theta(N)\cdot x}=
\bigcup_{x\Bro y} \theta(N)\cdot y$.
\item[{\rm (iv)}] $\Bro$ is a partial order.
\end{enumerate}
\end{theorem}

\begin{proof} (i)
Suppose that
$a^-\Smr x \Smr a^+$ and set
$S_\pm={\cal W}_{\cal M}^\pm(a^\mp)$.
Then the $\lambda$-Lemma implies $a^-\Smr a^+$.
This argument proves
that $\Smr$ is already transitive and hence
agrees with $\Smo$.

Next we show that
\[\overline{ {\cal W}_{\cal M}^+(x)}=
\bigcup_{y\preceq_{\cal M} x} {\cal W}_{\cal M}^+(y).\]
The inclusion $``\subseteq$'' follows since the right hand
side is an attractor for  $e^{-tX}$ and as such closed.
For the converse we again use the
$\lambda$-Lemma: Let $x\Smr y$ and
$z\in {\cal W}_{\cal M}^+(y)$. Choose a small
disc $S_-$ through $z$ intersecting
${\cal W}_{\cal M}^+(y)$ transversally.
Further we set  $S_+={\cal W}_{\cal M}^+(x)$. Then
$x \Smr y$ implies that $S_+$ intersects
${\cal W}_{\cal M}^-(y)$ (transversally), so that we find
a $t>0$ with $e^{tX}S_-\cap S_+\not=\emptyset$. But
then also
$S_-\cap S_+=S_-\cap e^{tX}S_+\not=\emptyset$. Since
$S_-$ was
arbitrarily small this proves the inclusion $``\supseteq$".

Now we see that $y\Smo x$ is equivalent to
$N\cdot y={\cal W}_{\cal M}^+(y)\subseteq
\overline{{\cal W}_{\cal M}^+(x)}=\overline{N\cdot x}$
which in turn just means $y\preceq_{Br} x$.
This proves (i) and (ii).

To prove (iii) we  observe that the time
reversed flow interchanges the roles of stable and unstable
manifolds and flips the Smale order around.
Therefore we only have to  apply  Theorem \ref{Bruhatcells1}.

 (iv) Suppose that $x\Bro y \Bro x$. Then
$x\Smr y\Smr x$ which implies
$x=y$ since  we have a gradient flow.

\end{proof}

We note that this theorem could also be derived from results of
Kazhdan and Lusztig (see \cite{De85}).

\section{Attractor-Repellor Pairs}\label{attrep}

In this section we determine  the attractor-repellor pairs for
special flows on complex flag manifolds and show how one can identify
the attractor network with the lattice of upper sets in ${\cal F}_{\cal M}$.
For details concerning attractor networks  we refer to
\cite{RobbinSalamon}.

\label{basicdefinitions}
We begin by reviewing some basic material
concerning attractors and repellors for general
continuous time dynamical systems.
An {\it attractor} for $\Phi$ in ${\cal M}$ is a compact
invariant set ${\cal A}\subseteq {\cal M}$ which admits a
neighborhood $U$ such that
\begin{equation}\label{attrnbhd1}
\Phi(t,\overline U)\subseteq {\rm int}\, U\quad\quad \forall t>0
\end{equation}
and
\begin{equation}\label{attrnbhd2}
{\cal A}=\bigcap_{t\ge 0}\Phi(t,U),
\end{equation}
where ${\rm int}\, U$ is the interior of $U$.
A {\it repellor} for $({\cal M},\Phi)$ is an attractor for the
time-reversed dynamical system.
We denote the set of all attractors for $\Phi$
by ${\cal A}({\cal M},\Phi)$ and the set of all
repellors for $\Phi$ by ${\cal A}({\cal M},\Phi^{-1})$.
Then ${\cal A}({\cal M},\Phi)$ is a lattice with respect to
the set theoretic meet and join operations.
For any subset ${\cal A}$ of ${\cal M}$ consider the set
$${\cal A}_{\fix}:=\{x\in {\cal A}\mid (\forall t\in \R) \Phi(t,x)=x\}$$
of fixed  points in ${\cal A}$.

\label{sinksource}

\begin{proposition} \label{fixofstable}
Assume that every flow line in ${\cal M}$ has a sink
and a source. In other words:
\begin{equation}\label{assumptionsinksource}
\lim_{t\to\pm\infty}\Phi(t,x)
\quad{\rm exists}\quad \forall x\in {\cal M}.
\end{equation}
Then for any closed invariant subset ${\cal A}\subseteq {\cal M}$
we have
\begin{enumerate}
\item[{\rm(i)}] ${\cal W}^{\pm}({\cal A})=
\bigcup_{x\in {\cal A}_{\fix}}{\cal W}^{\pm}(x)$.
\item[{\rm (ii)}] ${\cal W}^{\pm}({\cal A})_{\fix}={\cal A}_{\fix}$.
\end{enumerate}
\end{proposition}

 If ${\cal A}$ is  an attractor, then
${\cal A}^*={\cal M}\setminus {\cal W}^+({\cal A})$
is a repellor, called the {\it dual repellor} of ${\cal A}$.
Then $({\cal A},{\cal A}^*)$ is called an {\it attractor-repellor pair}.
From Proposition \ref{fixofstable} one derives

\begin{proposition} \label{attractorfix}
Let $({\cal A},{\cal A}^*)$ be an attractor-repellor pair.
Then we have
\begin{enumerate}
\item[{\rm (i)}] ${\cal A}^*_{\fix}=
{\cal M}_{\fix}\setminus {\cal A}_{\fix}$.
\item[{\rm(ii)}] ${\cal A}^*=
\bigcup_{x\in {\cal A}^*_{\fix}}{\cal W}^+(x)$ is closed.
\item[{\rm(iii)}] ${\cal A}=
\bigcup_{x\in {\cal A}_{\fix}}{\cal W}^-(x)$ is closed.
\item[{\rm(iv)}] ${\cal M}=
{\cal A}\cup\big({\cal W}^+({\cal A})\cap
{\cal W}^-({\cal A}^*)\big)\cup {\cal A}^*$
is a disjoint union.
\end{enumerate}
\end{proposition}

Recall from \cite{RobbinSalamon} that ${\cal A}({\cal M},\Phi)$ is a lattice
with respect to the usual set operations and that an {\it attractor network}
is a finite sublattice of ${\cal A}({\cal M},\Phi)$.
Before we can formulate and prove the general theorem that leads
the way to a description of
attractor networks for special flows in complex flag manifolds
we have to introduce some more order theoretic concepts.

To each partially ordered set $P$ one assigns the dual distributive lattice
$$
 P^* = {\rm Hom}_{\rm order}(P,\{0,1\})
$$
of order homomorphisms, where  the two point set $\{0,1\}$ is a
lattice with respect to  operations  $\max$ and $\min$.
A different way to view $P^*$ is to identify an element
$\alpha\in P^*$ with $\alpha^{-1}(1)$. This set is an upper
set with respect to the ordering. This means that for any
$p\in\alpha^{-1}(1)$ and $p\le q\in P$ we have
$q\in \alpha^{-1}(1)$. Conversely  each upper set
$U\subseteq P$ gives rise to an order homomorphism
$\chi_{U}\colon P\to \{0,1\}$, where
$\chi_{U}$ is the characteristic
function of $U$. In this way $P^*$ gets identified with
the lattice of upper subsets of $P$ under the usual set
theoretic join (union) and meet (intersection).

\begin{theorem}\label{attractorparam}
Suppose that the
 dynamical system $({\cal M},\Phi)$ is admissible with isolated fixed points
and $\Smo$ is a partial order.
Then all the attractor-repellor pairs $({\cal A},{\cal A}^*)$
are of the form
$({\cal A},{\cal A}^*)=({\cal A}_{\cal U},{\cal A}_{\cal U}^*)$
with
$${\cal A}_{\cal U}=\bigcup_{x\in {\cal U}}{\cal W}^-(x),\quad
{\cal A}_{\cal U}^*=\bigcup_{x\in {\cal L}}{\cal W}^+(x),$$
where ${\cal U}$ is a compact isolated upper subset of ${\cal M}_{\fix}$
with respect to the Smale order
and ${\cal L}={\cal M}_{\fix}\setminus {\cal U}$.
The map
$${\cal M}_{\fix}^*\mapsto {\cal A}({\cal M},\Phi), \quad {\cal U}\mapsto{\cal A}_{\cal U}
$$
is a lattice isomorphism.
\end{theorem}

\begin{proof} Given an attractor ${\cal A}$ and an attracting neighbourhood
$U\subseteq {\cal M}$ as described in (\ref{attrnbhd1}) and
 (\ref{attrnbhd2}) we see that
${\cal L}:=U\cap{\cal M}_{\fix}= {\cal A}_{\fix}$.
This implies that ${\cal L}$ is an open and closed lower
set in ${\cal M}_{\fix}$. We set ${\cal U}={\cal M}_{\fix}\setminus{\cal L}$.
Then Proposition \ref{attractorfix} proves the equalities
${\cal A}_{\cal U}={\cal A}$ and ${\cal A}_{\cal U}^*={\cal A}^*$.
Conversely, since there are only finitely many open and closed
subsets in ${\cal M}_{\fix}$ we can apply \cite[Thm.~1.3]{RobbinSalamon}
to see that each attractor has to be of this form.
The last claim is now clear.
\end{proof}

 \label{Lyapunovmap}
Theorem \ref{attractorparam} says that for the
partially ordered  (with respect to the Smale order) set
${\cal M}_{\fix}$ of  isolated  fixed points of  $\Phi(t,x)$
the dual lattice is the lattice of attractors of the flow.
The correspondence is given by
$$
     {\cal A}_\alpha =
\bigcup_{\alpha(p)=0} {\cal W}^-(p)
$$
for $\alpha\in {\cal M}_{\fix}^*$.
For each attractor the dual repellor is given by
$$
      {\cal A}_\alpha^* =
\bigcup_{\alpha(p)=1} {\cal W}^+(p).
$$

\section{The Attractor  Network for Special Flows on Complex Flag Manifolds}

\label{specialattractors}

Let $V$ be finite dimensional complex vector space with an inner
product $\la\cdot\,|\cdot \ra$ and $\varphi\in \End(V)$ a
selfadjoint linear map. For each eigenvalue $\nu\in \spec(\varphi)$
of $\varphi$ we set $V_+^{(\nu)}:=\sum_{\nu\le \mu}V_\mu$ and
$V_-^{(\nu)}:=\sum_{\mu<\nu}V_\mu.$ Then $V=V_+^{(\nu)}\oplus
V_-^{(\nu)}$.

\begin{proposition} The attractor-repellor pairs for
$(\Proj(V),\Phi)$ are given by
\[{\cal A}_\nu:=\Proj\left(V_+^{(\nu)}\right)
\quad\mbox{and}\quad {\cal A}^*_\nu:=\Proj\left(V_-^{(\nu)}\right)\]
for $\nu\in \spec(\varphi)$.
\end{proposition}

\begin{proof} This follows from (\ref{tildeM1}), Remark \ref{tildeM2}
and Theorem \ref{attractorparam}.
\end{proof}

We return to the particular case of a complex flag manifold ${\cal
M}=G\cdot x_c\subseteq \Proj(V)$ and a special flow $\Phi_X$ on
${\cal M}$ considered in \S~\ref{specialflows}.   Theorem
\ref{Bruhatorder1} shows that the hypotheses for Theorem
\ref{attractorparam} are satisfied and that the Smale order agrees
with the Bruhat order.

\label{attractor}
Let ${\cal A}$ be an attractor in ${\cal M}$ for the flow
$\Phi_X$ and ${\cal A}^*$ its dual repellor. Then, by  Proposition  \ref{fixofstable}, we have
$${\cal A}={\cal W}^{-}({\cal A})=
\bigcup_{x\in {\cal A}_{\fix}}{\cal W}^{-}(x)=
\bigcup_{x\in {\cal A}_{\fix}}(\theta N)\cdot x$$
and
$${\cal A}^*={\cal W}^{+}({\cal A}^*)=
\bigcup_{x\in {\cal A}_{\fix}^*}{\cal W}^{+}(x)=
\bigcup_{x\in {\cal A}_{\fix}^*}N\cdot x.$$
Since  repellors are closed, the Theorem \ref{Bruhatorder1} shows that
${\cal A}_{\fix}^*\subseteq {\cal M}_{\fix}$ is a {\it lower set} with
respect to the Bruhat
order.
 This means that if $x\in {\cal A}_{\fix}^*$,
then all $y\in {\cal M}_{\fix}$ with $y\Bro x$
also belong to ${\cal A}_{\fix}^*$.
Similarly
${\cal A}_{\fix}\subseteq {\cal M}_{\fix}$ is an {\it upper set} with respect to the Bruhat
order. In fact, it turns out that {\it all} lower and upper
sets can be obtained in this way.

\begin{theorem} \label{attractor1}
Let ${\cal L}({\cal F}_{\cal M})$ be the lattice
of lower sets in ${\cal F}_{\cal M}$ and $\ {\cal U}({\cal F}_{\cal M})$
be the lattice of upper sets in ${\cal F}_{\cal M}$. Then the maps
\[{\cal U}({\cal F}_{\cal M})\to{\cal A}({\cal M},\Phi),\quad
U\mapsto {\cal A}_U:=\bigcup_{x\in U} (\theta N)\cdot x\]
and
\[{\cal L}({\cal F}_{\cal M})\to{\cal A}({\cal M},\Phi^{-1}),\quad
L\mapsto {\cal R}_{L}:=
\bigcup_{x\in L} N\cdot x\]
are lattice isomorphisms and $ {\cal R}_{L}$ is the dual
repellor of ${\cal A}_{{\cal F}_{\cal M}\setminus L}$.
\end{theorem}

\begin{proof}
In view of Theorem \ref{attractorparam} only the last claim remains to be proved.
To do that, let $L={\cal F}_{\cal M}\setminus U$. Then from
Proposition  \ref{attractorfix} we see
\begin{eqnarray*}
{\cal M}&=&\bigcup_{(x,x')\in {\cal F}_{\cal M}\times
{\cal F}_{\cal M}} N\cdot x\cap (\theta N)\cdot x',\\
{\cal A}_U&=&\bigcup_{(x,x')\in {\cal F}_{\cal M}\times
U} N\cdot x\cap (\theta N)\cdot x',\\
{\cal R}_L&=&\bigcup_{(x,x')\in L\times
{\cal F}_{\cal M}} N\cdot x\cap (\theta N)\cdot x',\\
{\cal W}^s({\cal A}_U)\cap {\cal W}^u({\cal R}_L)&=&
\bigcup_{(x,x')\in U\times L} N\cdot x\cap (\theta N)\cdot x',
\end{eqnarray*}
so that we have a disjoint union
$${\cal M}={\cal A}_U\cup
\big({\cal W}^s({\cal A}_U)\cap {\cal W}^u({\cal R}_L)\big)
\cup {\cal R}_L.$$
Now \cite[Prop.\ 1.4]{RobbinSalamon} shows that
${\cal A}_U$ is an attractor and ${\cal R}_L={\cal A}_U^*$
is its dual repellor.

\end{proof}

\end{document}